\documentclass{elsart}
\newcounter{guy}
\newtheorem{theorem}{Theorem}[section]
\newtheorem{lemma}{Lemma}[section]
\newtheorem{corollary}{Corollary}[section]
\newtheorem{proposition}{Proposition}[section]
\usepackage{amsfonts}
\usepackage{amsmath}
\newcommand{\beats}{\rightarrow}
\begin{document}
    
    \begin{frontmatter}
    
    \title{Quadrangularity in Tournaments}
    \author[UCD]{J. Richard Lundgren}, 
    \author[Bris]{Simone Severini} and
    \author[UCD]{Dustin J. Stewart\corauthref{cor1}}
    \ead{dstewart@math.cudenver.edu}
    \corauth[cor1]{Corresponding author}
    
    \address[UCD]{University of Colorado at Denver, Denver, 
    CO 80217}
    \address[Bris]{University of Bristol, Bristol, United Kingdom}

    \begin{abstract} 
    The pattern of a matrix $M$ is a (0,1)-matrix which replaces all 
    non-zero entries of $M$ with a 1.  There are several contexts 
    in which studying the patterns of orthogonal matrices can be 
    useful.  One necessary 
    condition for a matrix to be orthogonal is a property known as 
    combinatorial orthogonality.  If the adjacency matrix of a
    directed graph forms a pattern of a combinatorially orthogonal 
    matrix, we say the digraph is quadrangular.  We look at 
    the quadrangular property in tournaments and regular tournaments. 
    \end{abstract}
    
    \begin{keyword} combinatorial orthogonality, quadrangular, 
    tournament \end{keyword}
    
    \end{frontmatter}

    \section{Introduction}
    
    A \em directed graph \em or \em digraph, \em  $D,$ 
    is a set of vertices, $V(D)$, 
    together with a set of ordered pairs of the vertices, $A(D)$, 
    called arcs.  If $(u,v)$ is an arc in a digraph, we say 
    that $u$ \em beats \em $v$ or $u$ \em dominates \em $v$, and typically 
    write this as $u\beats v$.
    If $v\in V(D)$ then we define the \em outset \em of $v$ by,
    $$O_{D}(v)=\{u\in V(D):(v,u)\in A(D)\}.$$
    That is, $O_{D}(v)$ is all vertices in $D$ which $v$ beats.  
    Similarly, we 
    define the set of all vertices in $D$ which beat $v$ to be the \em 
    inset \em 
    of $v$, written,
    $$I_{D}(v)=\{u\in V(D):(u,v)\in A(D)\}.$$
    The \em in-degree \em and \em out-degree \em of a 
    vertex $v$ are $d_{D}^{-}(v)=|I_{D}(v)|$ and $d_{D}^{+}(v)=|O_{D}(v)|$ 
    respectively.  When it is clear to which digraph $v$ belongs, we 
    will drop the subscript.
    The minimum out-degree of $D$ is the smallest 
    out-degree of any vertex in $D$, and is denoted by $\delta^{+}(D).$  
    We define the minimum in-degree 
    of $D$ similarly, and represent it by 
    $\delta^{-}(D)$. 
    A \em tournament, \em $T$, is directed graph with the property 
    that for each $u\neq v\in V(T)$ exactly one of $(u,v)$, $(v,u)$ is 
    in $A(T)$, and $(u,u)\not\in A(T)$.

    \medskip\indent
    Let $x=(x_{1},x_{2},\ldots,x_{n})$ and $y=(y_{1},y_{2},\ldots,y_{n})$ 
    be $n-$vectors over any field.  We say that $x$ and 
    $y$ are \em combinatorially orthogonal \em if $|\{i:x_{i}y_{i}\neq 
    0\}|\neq 1.$  Let $M$ be an $n\times n$ matrix.  If every two rows 
    of $M$ are combinatorially orthogonal, then we say that $M$ is \em 
    combinatorially row-orthogonal. \em  If both $M$ and $M^{T}$ have 
    the property of being combinatorially row-orthogonal, then we say 
    that $M$ is a \em combinatorially orthogonal \em matrix.  
    In \cite{Klee}, Beasley, Brualdi and Shader study the 
    combinatorial orthogonality property to determine that an 
    orthogonal matrix which cannot be decomposed into two smaller 
    orthogonal matrices must contain at least $4n-4$ non-zero 
    entries.  They also determine a family of matrices so that any 
    combinatorially orthogonal matrix which cannot be decomposed and 
    meets this bound belongs to 
    this family, up to arbitrary row and column permutations.

    \medskip\indent
    Given an $n\times n$ matrix $M$, the \em pattern \em of $M$ is 
    a $n \times n$ matrix $M',$ defined by
    $$M'_{i,j}=\left\{\begin{array}{rl} 0 & \mbox{ if $M_{i,j}=0$} \\
    1 & \mbox{ if $M_{i,j}\neq 0$}\end{array}\right.$$
    where $M_{i,j}$ denotes the $i,j$ entry of $M.$
    It follows quickly from the definition that a matrix is 
    combinatorially orthogonal if and only if the associated pattern 
    matrix is combinatorially orthogonal.  
    
    \medskip\indent
    Using the pattern of a matrix allows us to relate the concept of 
    combinatorial orthogonality to digraphs.
    We define the \em adjacency matrix \em $A$ of a digraph $D$, 
    with $V(D)=\{v_{1},v_{2},\ldots,v_{n}\}$ by,
    $$A_{i,j}=\left\{\begin{array}{rl}1 & \mbox{ if $v_{i}\beats v_{j}$, 
    and}\\ 
    0 & \mbox{ otherwise.} \end{array}\right.$$
    If the adjacency matrix of $D$ is the pattern of a matrix $M$, we 
    simplify the statement by saying that $D$ is the digraph of $M$.
    
    \medskip\indent
    Note, $D$ is the digraph of a combinatorially 
    orthogonal matrix if and only if for all $u\neq v\in V(D),$ 
    $|O(u)\cap O(v)|\neq 1$ and $|I(u)\cap I(v)|\neq 1.$  
    We call a digraph with these properties \em quadrangular. \em   
    If we only require $|O(u)\cap O(v)|\neq 1$ for all 
    $u\neq v\in V(D),$ we say $D$ is \em out-quadrangular. \em   
    Similarly, 
    if $|I(u)\cap I(v)|\neq 1$ for all $u\neq v\in V(D),$ we say that $D$ 
    is \em in-quadrangular. \em  
    In \cite{Gibson}, Gibson and Zhang study a similar quadrangular 
    property 
    by looking at combinatorial orthogonality in the reduced 
    adjacency matrices of bipartite graphs.  By citing a theorem of Reid 
    and Thomassen \cite{Reid/Thomassen}, Gibson and Zhang give a graph 
    theoretic proof of the bound of Beasley, Brualdi and Shader found in 
    \cite{Klee}. 
    In this paper we are interested in tournaments which 
    have the quadrangular property.

    \medskip\indent
    Characterizing digraphs of orthogonal matrices is a method to unveil 
    their
    combinatorial properties. An understanding of combinatorial properties
    of orthogonal matrices could be useful in approaching the existence 
    problem
    for weighing matrices.
    Also, it may provide insight in contexts where combinatorial
    objects and orthogonal matrices
    naturally appear. For example, in the theory of quantum
    computation and information (\cite{Kus}).

   \section{Quadrangular tournaments}
   
   \medskip\indent
   In this section we give some classifications of quadrangular 
   tournaments, but 
   first we need a few more definitions.
   Let $T$ be a tournament.  We obtain the \em dual of $T,$ \em 
   which we represent by $T^{r},$ by 
   forming the tournament on the same vertices of $T$ with $x\beats y$ 
   in $T^{r}$ if and only if $y\beats x$ in $T$.  We use the 
   notation $T^{r}$ because this is sometimes referred to as the 
   reversal of $T$.  Also, a \em transmitter \em in $T$ is a vertex 
   which dominates all other vertices of $T$, and a \em receiver \em is a 
   vertex which is dominated by all other vertices of $T$.
   
   \medskip\indent
   A \em dominant pair \em in a tournament $T$ is a pair of vertices 
   $u,v$ so that every other vertex in $T$ is dominated by at least one of 
   $u$ or $v$.  The \em domination 
   graph \em of $T$, denoted by dom$(T)$ is the graph formed on the 
   same vertices of $T$ with an edge between $x,y$ if and only if $x$ 
   and $y$ form a dominant pair in $T.$  The \em competition graph \em 
   of $T$ is the graph formed on the same vertices of $T$ with an edge 
   between $x$ and $y$ if and only if there exists some $z$ such that 
   $x\beats z$ and $y\beats z.$  Fisher, Lundgren, Merz and Reid 
   \cite{Lundgren}
   showed that the domination graph of $T$ is isomorphic to the 
   competition graph of $T^{r}.$

   \medskip\indent
   A \em dominating set \em in a digraph $D$ is a set of vertices $S$ 
   such that every vertex in $D$ is in $S$ or dominated by a vertex 
   in $S$.  The \em domination number \em of a 
   digraph, $\gamma(D),$ is the size of a smallest dominating set 
   in $D.$ 
   Note that a dominating set of size $2$ in a 
   tournament is a dominant pair.  So, if $\gamma(T)>2$ then $T$ has 
   no dominating pairs and so $E($dom$(T))=\emptyset.$  We now use 
   these concepts to classify some quadrangular tournaments.

   \begin{theorem}\label{tourn2}
   Let $T$ be a tournament on $3$ or more vertices 
   with a transmitter $s$ and  
   receiver $t$.  Then $T$ is quadrangular if and only if both 
   $\gamma(T-\{s,t\})>2$ and 
   $\gamma((T-\{s,t\})^{r})>2.$  \end{theorem} 
   
   \begin{pf*}{Proof.} Let $T$ be a tournament with a transmitter $s$ and 
   receiver 
   $t.$  Suppose that both $\gamma(T-\{s,t\})>2$ and 
   $\gamma((T-\{s,t\})^{r})>2.$  Then, 
   $E($dom$(T-\{s,t\}))=E($dom$((T-\{s,t\})^{r}))=\emptyset.$  Thus the 
   competition graphs of both $T-\{s,t\}$ and $(T-\{s,t\})^{r}$ are 
   complete.  That is, for all $x,y\in V(T-\{s,t\})$ there exist 
   $w,z\in V(T)$ such that $w\beats x,$ $w\beats y,$ $x\beats z$ 
   and $y\beats z.$  
   Pick $u\neq v\in V(T).$  We consider three cases.
   
   \medskip\noindent
   \em Case 1: \em Suppose  $u,v\not\in\{s,t\}$. 
   Then, as noted before, 
   there exist vertices $w,z\in V(T-\{s,t\})$ so that 
   $z\in O(u)\cap O(v)$ and $w\in I(u)\cap I(v).$  Also, $s\in 
   O(u)\cap O(v)$ and $t\in I(u)\cap I(v).$  So, $|O(u)\cap O(v)|\geq 
   2$ and $|I(u)\cap I(v)|\geq 2.$  
   
   \medskip\noindent
   \em Case 2: \em  Now assume that one of $u$ or $v$ is $t,$ say $u=t.$  
   Since $O(t)=\emptyset,$ $O(t)\cap O(v)=\emptyset,$ so $|O(t)\cap 
   O(v)|=0.$  Also, $I(t)=V(T)-t,$ so $I(t)\cap I(v)=I(v).$  If 
   $v=s,$ then $I(v)=\emptyset,$ thus $|I(t)\cap I(v)|=0.$  So, suppose 
   $v\neq s.$  Since $\gamma(T-\{s,t\})>2,$ there exists $w\in 
   V(T-\{s,t\})$ such that $w\beats v,$ for otherwise $v$ would be a 
   dominating set of size $1$ in $T-\{s,t\}.$  Thus, 
   $s,w \in I(v),$ and $|I(t)\cap I(v)|\geq 2$ as desired.     
   
   \medskip\noindent
   \em Case 3: \em Now, assume that one of $u,v$ is $s,$ say $u=s.$  
   Since $I(s)=\emptyset,$ $I(s)\cap I(v)=\emptyset,$ so $|I(s)\cap 
   I(v)|=0.$  Also, since $O(s)=T-s,$ $O(s)\cap O(v)=O(v).$  The case 
   with $v=t$ was covered in case 2, so assume $v\neq t.$  Since 
   $\gamma((T-\{s,t\})^{r})>2$ there exists $w$ such that $v\beats w,$ for 
   otherwise $w$ would form a dominating set of size $1$ in 
   $(T-\{s,t\})^{r}.$  So, 
   $w,t \in O(v),$ and so $|O(s)\cap O(v)|\geq 2.$
   
   \medskip\indent
   Now assume that $T$ is a quadrangular tournament with both a 
   transmitter $s$ and receiver $t.$  If $u,v\in V(T-\{s,t\}),$ then 
   $s\in O(u)\cap O(v),$ and $t\in I(u)\cap 
   I(v).$  Since $T$ is quadrangular and $|O(u)\cap O(v)|\geq 1$ 
   and $|I(u)\cap I(v)|\geq 1,$ there must exist vertices $w,z$ in 
   $T-\{s,t\}$ such that $z\in O(u)\cap O(v)$ and $w\in I(u)\cap 
   I(v).$  Since $w$ beats $u$ and $v$ they cannot be a dominant pair 
   in $T-\{s,t\}$ and since $u$ and $v$ beat $z$ they cannot be a 
   dominant pair in $T-\{s,t\}.$  Thus, 
   $E($dom$(T-\{s,t\}))=E($dom$((T-\{s,t\})^{r})) =\emptyset.$  
   Equivalently, $\gamma(T-\{s,t\})>2$ and 
   $\gamma((T-\{s,t\})^{r})>2.$  This completes the proof.\qed \end{pf*}

   \medskip\indent
   It was shown by Fisher et.al. in \cite{Fisher}
   that a tournament on fewer than $7$ 
   vertices must contain a dominant pair.  It is known that 
   the quadratic residue tournament on $7$ vertices, $QR_{7}$, has 
   domination number $3$, and $QR_{7}$ is isomorphic to its dual, so 
   $\gamma(QR_{7}^{r})=3$.  Thus a tournament $T$  
   on $9$ vertices
   with a transmitter $s$ and receiver $t$ such that 
   $T-\{s,t\}=QR_{7}$ is the smallest example of a quadrangular 
   tournament with both a transmitter and receiver.  We now consider 
   the case when a tournament has a transmitter or receiver, but not 
   both.

   \begin{theorem}\label{tourn3} Let $T$ be a tournament with a 
   transmitter $s$ and no receiver.  Then $T$ is quadrangular if and only 
   if, 
   $\gamma(T-s)>2$, $T-s$ is out-quadrangular, and $\delta^{+}(T-s)\geq 
   2.$  \end{theorem}
   
   \begin{pf*}{Proof.} First suppose that $\gamma(T-s)>2,$ $T-s$ is 
   out-quadrangular, and $\delta^{+}(T-s)\geq 2.$    Pick $u\neq v\in 
   V(T).$  
   First suppose that $u,v\in T-s.$
   Since $\gamma(T-s)>2$ there exists a $x\in V(T-s)$ such that 
   $x\beats u$ and $x\beats v.$  So, $s,x\in I(u)\cap I(v)$ 
   and so $|I(u)\cap I(v)|\geq 2.$  Also, since $T-s$ is 
   out-quadrangular $|O(u)\cap O(v)|\neq 1.$  Now, suppose that one 
   of $u,v$ is $s,$ say $u=s.$  Since $I(s)=\emptyset,$ $|I(s)\cap 
   I(v)|=0.$  Also, since $\delta^{+}(T-s)\geq 2,$ 
   $|O(s)\cap O(v)|=|O(v)|\geq 2,$
   Thus, $T-s$ is quadrangular as desired.
   
   \medskip\indent
   Now, assume that $T$ is quadrangular.  Since $O(s)=T-s,$ $|I(u)\cap 
   I(v)|\geq 1$ for all $u,v\in V(T-s).$  Since $T$ is quadrangular 
   this means we must have $|I(u)\cap I(v)|\geq 2$ for each 
   $u\neq v\in V(T-s).$  Thus, for all $u\neq v\in V(T-s)$,
   there must exist some $x\in V(T-s)$ such that $x\beats u$ and 
   $x\beats v.$  So, $\gamma(T-s)>2.$  Since $T$ has no receiver, 
   $|O(v)|\geq 1$ for all $v\in V(T).$  Since 
   $O(s)\cap O(v)=O(v)$ for all $v\in V(T-s),$ and $T$ is 
   quadrangular, we must then have that 
   $$|O(v)|=|O(s)\cap O(v)|\geq 2.$$
   Thus, $\delta^{+}(T-s)\geq 2.$  Now, pick $u\neq v\in V(T-s).$  
   Since $T$ is quadrangular, $|O(u)\cap O(v)|\neq 1$.  So $T-s$ is 
   out-quadrangular.  
   \qed \end{pf*}

   \medskip\indent
   If $T$ is a tournament 
   with a receiver and no transmitter, then it is 
   the dual of a tournament with a transmitter and no receiver.  
   Obviously, 
   a tournament is quadrangular if and only if its dual is.  So, 
   by Theorem~\ref{tourn3}, $T$ is quadrangular if and only if 
   $\gamma((T-t)^{r})>2$, $(T-t)^{r}$ is out-quadrangular and 
   $\delta^{+}((T-t)^{r})\geq 2.$  Since $(T-t)^{r}$ being 
   out-quadrangular is equivalent to $T-t$ being in-quadrangular, 
   and  
   $\delta^{+}((T-t)^{r})=\delta^{-}(T-t)$ we get the following 
   corollary.

   \begin{corollary}\label{tourn4} Let $T$ be a tournament with a 
   receiver $t$
   and no transmitter.  Then $T$ is quadrangular if and only if 
   $\gamma((T-t)^{r})>2$, $T-t$ is in-quadrangular, and 
   $\delta^{-}(T-t)\geq 2.$ \end{corollary}

   \medskip\indent
   A tournament is called \em strongly connected \em if any two 
   vertices in the tournament are mutually reachable by a directed 
   path.  If a tournament is not strongly connected, then it can be 
   partitioned into maximal strongly connected components.  
   Further, these
   strong components can be labeled $T_{1},T_{2},\ldots,T_{m}$ such that 
   every vertex in $T_{i}$ beats every vertex of $T_{j}$ whenever 
   $i<j.$  The component $T_{1}$ is called the \em initial strong 
   component \em and the component $T_{m}$ is called the \em terminal 
   strong component. \em

   \begin{theorem}\label{tourn5} Let $T$ be a tournament with no 
   transmitter or receiver which is not strongly connected.  Then $T$ is 
   quadrangular if and only if the initial strong component, $T_{1}$, 
   is in-quadrangular with $\delta^{-}(T_{1})\geq 2$ and the terminal 
   strong component, $T_{m}$, is out-quadrangular with 
   $\delta^{+}(T_{m})\geq 2.$  \end{theorem}
   
   \begin{pf*}{Proof.} Let $T$ be a tournament with no 
   transmitter or receiver, which is 
   not strongly connected.  Suppose 
   that $T_{1}$ is in-quadrangular with $\delta^{-}(T_{1})\geq 2,$ 
   and that $T_{m}$ is out-quadrangular with $\delta^{+}(T_{m})\geq 
   2.$  Note also that since $T$ has no transmitter or receiver, 
   $T_{1}$ and $T_{m}$ must contain at least $3$ vertices each.
   Pick $u\neq v\in V(T).$  We consider $5$ cases.
   
   \medskip\noindent
   \em Case 1: \em Suppose that $u$ and $v$ are in neither 
   $T_{1}$ nor $T_{m}.$  Every vertex of $T_{1}$ beats 
   every vertex in $T-T_{1}$ and every vertex of $T_{m}$ is beaten by 
   every vertex of $T-T_{m}.$  So, since $T$ has no transmitter or 
   receiver,
   $$|O(u)\cap O(v)|\geq |V(T_{m})|\geq 3 \mbox{ and } |I(u)\cap 
   I(v)|\geq |V(T_{1})|\geq 3.$$
   
   \medskip\noindent
   \em Case 2: \em Suppose that both $u,v\in T_{1}.$  Then, since 
   $T_{1}$ is in-quadrangular, $|I(u)\cap I(v)|\neq 1.$  Also, $u$ 
   and $v$ beat every vertex in $T-T_{1}$, in particular, 
   $T_{m}\subseteq O(u)\cap O(v).$  Thus, since $T$ has no receiver,
   $$|O(u)\cap O(v)|\geq |V(T_{m})|\geq 3.$$
   
   \medskip\noindent
   \em Case 3: \em Suppose that both $u,v\in T_{m}.$  Then, since 
   $T_{m}$ is out-quadrangular, $|O(u)\cap O(v)|\neq 1.$  Also, since 
   $T$ has no transmitter, $|I(u)\cap I(v)|\geq |V(T_{1})|\geq 3.$
   
   \medskip\noindent
   \em Case 4: \em Suppose that $u\in T_{1}$ and $v\not\in T_{1}.$  
   Since $v\not\in T_{1}$ we know that $I(u)\subseteq I(v)$ and so 
   $I(u)\cap I(v)=I(u).$  So, since $\delta^{-}(T_{1})\geq 2$, 
   $|I(u)\cap I(v)|=|I(u)|\geq 2$.  
   Also, since $u\in T_{1}$ and $v\not\in T_{1},$ we know that 
   $O(v)\subseteq 
   O(u)$.  Thus, $O(u)\cap O(v)=O(v)$.  If $v\not\in T_{m}$, then 
   $T_{m}\subseteq O(v),$ and so $|O(u)\cap O(v)|\geq |T_{m}|\geq 
   3$.  So, assume that $v\in T_{m}.$  Then $|O(u)\cap 
   O(v)|=|O(v)|\geq \delta^{+}(T_{m})\geq 2.$
   
   \medskip\noindent
   \em Case 5: \em Suppose that $u\in V(T_{m})$ and $v\not\in V(T_{m})$.  
   Since $v\not\in V(T_{m})$, $O(u)\subseteq V(T_{m})\subseteq O(v),$ and 
   so $O(u)\cap O(v)=O(u).$  So, since $\delta^{+}(T_{m})\geq 2,$ 
   $|O(u)\cap O(v)|=|O(v)|\geq 2.$  Now, if $v\in V(T_{1})$ then we 
   showed in case 4 that $|I(u)\cap I(v)|\geq 2.$  So, assume that 
   $v\not\in V(T_{1})$.  Then, every vertex in $T_{1}$ beats both $u$ and 
   $v$, and so $|I(u)\cap I(v)|\geq |V(T_{1})|\geq 3$.

   \medskip\indent
   Now, assume that $T$ is quadrangular.  Since $T$ is quadrangular, 
   if $u\neq v\in V(T_{1})$, then $|I(u)\cap I(v)|\neq 1,$ so $T_{1}$ is 
   in-quadrangular.  Also, since $T$ is quadrangular, 
   if $u\neq v\in V(T_{m})$ 
   then $|O(u)\cap O(v)|\neq 1,$ and so $T_{m}$ is out-quadrangular.   
   Now, pick $u\in V(T_{1})$ and $v\in V(T_{m})$.  Then, $O(u)\cap 
   O(v)=O(v),$ and $I(u)\cap I(v)=I(u).$  
   Since $T$ has no receiver, $|O(v)|\geq 1$, and so we must 
   have that $|O(v)|=|O(u)\cap O(v)|\geq 2$.  Thus, 
   $\delta_{+}(T_{m})\geq 2$.  Also, since $T$ has no transmitter, 
   $|I(u)|\geq 1,$ and so $|I(u)|=|I(u)\cap I(v)|\geq 2$.
   Thus, 
   $\delta^{-}(T_{1})\geq 2$.  These are the conditions from 
   the theorem statement, and so the result follows.  
   \qed \end{pf*}
   
   \medskip\indent 
   We now give a characterization of quadrangular tournaments with 
   minimum in-degree $1$ or minimum out-degree $1.$ 
   First we need some lemmas.
   
   \begin{lemma}\label{tourn6} Let $T$ be a quadrangular tournament 
   with a vertex $x$ of out-degree $1$.  Say $x\beats y$ then 
   $O(y)=T-\{x,y\}$ \end{lemma}
   
   \begin{pf*}{Proof.} Suppose there exists a vertex 
   $v$ in $T-\{x,y\}$ such that 
   $v\beats y.$  Then, since $O(x)=y,$ $|O(x)\cap O(v)|=|\{y\}|=1.$  This 
   contradicts quadrangularity of $T.$  \qed \end{pf*}
   
   \medskip\indent
   Applying Lemma~\ref{tourn6} to the dual of $T$ we obtain the 
   following lemma.
   
   \begin{lemma}\label{tourn6a}  Let $T$ be a quadrangular tournament 
   with a vertex $x$ of in-degree $1$.  Say $y\beats x$, then 
   $I(y)=T-\{x,y\}.$  \end{lemma}

   \begin{theorem}\label{tourn7}  Let $T$ be a tournament on $4$ or 
   more vertices with  
   a vertex $x$ of out-degree $1,$ and say $x\beats y$.
   Then, $T$ is quadrangular if and only if 
   \begin{list}{\arabic{guy}.}{\usecounter{guy}}
   \item $O(y)=T-\{x,y\},$ 
   \item $\gamma(T-\{x,y\})>2$, 
   \item $\gamma((T-\{x,y\})^{r})>2,$ 
   \item $\delta^{+}(T-\{x,y\})\geq 2,$ 
   \item $\delta^{-}(T-\{x,y\})\geq 2.$  
   \end{list}\end{theorem}
   
   \begin{pf*}{Proof.} First, suppose that $T$ is quadrangular.  Then, by 
   Lemma~\ref{tourn6}, $O(y)=T-\{x,y\}.$  Now, pick vertices $u\neq v$ 
   in $T-\{x,y\}$.  Since $x\in O(u)\cap O(v)$ there must exist some 
   other vertex $w$ in $T-x$ for which $w\in O(u)\cap O(v).$  Since 
   $O(y)=T-\{x,y\},$ this $w$ must be in $T-\{x,y\}.$  So, there 
   exits $w\in T-\{x,y\}$ such that $w\in O(u)\cap O(v).$  This is 
   equivalent to saying $\gamma((T-\{x,y\})^{r})>2.$  Also, $y\in 
   I(u)\cap I(v).$  So, since $T$ is quadrangular, there must exist a 
   vertex $z$ in $T-y$ such that $z\in I(u)\cap I(v).$  Since 
   $O(x)=y,$ this vertex must be in $T-\{x,y\}.$  So, we must also 
   have that $\gamma(T-\{x,y\})>2.$  
   Now, if $v\in V(T-\{x,y\})$, 
   then $I(v)\cap I(x)= I(v).$  Since $\gamma(T-\{x,y\})>2,$ 
   $I(v)-y\neq \emptyset,$  So, $|I(v)|=|I(v)\cap I(x)|\geq 2.$  
   Thus, $\delta^{-}(T-\{x,y\})\geq 2.$  Also, $O(v)\cap O(y)=O(v)-x.$ 
   Since $\gamma((T-\{x,y\})^{r})>2,$ we have $O(v)-x\neq \emptyset,$ and 
   so 
   $|O(v)\cap O(y)|\geq 1.$  Thus, since $T$ is quadrangular, 
   $|O(v)|=|O(y)\cap O(v)|\geq 2.$  Thus, $\delta^{+}(T-\{x,y\})\geq 
   2.$  So, these conditions are necessary.
   
   \medskip\indent
   Now assume that $T$ is a tournament with a vertex $x$ such that 
   $O(x)=y,$ 
   and $O(y)=T-\{x,y\},$ $\gamma(T-\{x,y\})>2$, 
   $\gamma((T-\{x,y\})^{r})>2,$ $\delta^{+}(T-\{x,y\})\geq 2$ and 
   $\delta^{-}(T-\{x,y\})\geq 2.$  Pick $u\neq v\in V(T).$  We will  
   show $T$ is quadrangular using three cases.
   
   \medskip\noindent
   \em Case 1: \em Suppose $u,v\in V(T-\{x,y\}).$  Then, $x\in 
   O(u)\cap O(v),$ and since $\gamma((T-\{x,y\})^{r})>2,$ there exits 
   $w\in V(T-\{x,y\})$ such that $w\in O(u)\cap O(v).$  Thus, 
   $|O(u)\cap O(v)|>1.$  Also, $y\in I(u)\cap I(v),$ and since 
   $\gamma(T-\{x,y\})>2$ there exists $z\in V(T-\{x,y\})$ such that 
   $z\in I(u)\cap I(v).$  So, $|I(u)\cap I(v)|>1.$
   
   \medskip\indent
   \em Case 2: \em Suppose that $u=x.$  Then $O(u)=y$ and since 
   $y\not\in O(v),$ $|O(u)\cap O(v)|=0.$  Now, $I(u)\cap I(v)=I(v)-y$ 
   since $u=x.$  So, since $\delta^{-}(T-\{x,y\})\geq 2,$ $|I(u)\cap 
   I(v)|=|I(v)-y|\geq 2.$ 
   
   \medskip\indent
   \em Case 3: \em Suppose that $u=y.$  Then, $I(u)=x$ and since 
   $x\not\in I(v),$ $|I(u)\cap I(v)|=0.$  Now, $O(u)\cap O(v) 
   =O(v)-x$ since $u=y.$  So, since $\delta^{+}(T-\{x,y\})\geq 2,$ 
   $|O(u)\cap O(v)|=|O(v)-x|\geq 2.$  Thus, $T$ is quadrangular.
   \qed \end{pf*}
   
   \medskip\indent
   Applying Theorem~\ref{tourn7} to the dual of $T$ we obtain the 
   following corollary.
   
   \begin{corollary}\label{tourn7a} Let $T$ be a tournament with a 
   vertex $x$ with in-degree $1.$  Let $y=I(x).$  Then, $T$ is  
   quadrangular if and only if $I(y)=T-\{x,y\},$ 
   $\gamma(T-\{x,y\})>2$, $\gamma((T-\{x,y\})^{r})>2$, 
   $\delta^{+}(T-\{x,y\})\geq 2,$ and $\delta^{-}(T-\{x,y\})\geq 2.$  
   \end{corollary}

   \medskip\indent
   We now consider tournaments whose minimum out-degree 
   and in-degree are at least $2$.

   \begin{theorem}\label{odeg23} Let $T$ be an out-quadrangular 
   tournament and choose $v\in V(T)$.  Let $W$ be the sub-tournament 
   of $T$ induced on the vertices of $O(v)$.  Then $W$ contains no 
   vertices of out-degree $1$.  \end{theorem}
   
   \begin{pf*}{Proof.} Let $T$ be an out-quadrangular tournament, 
   and choose a vertex $v\in V(T)$.
   Let $W$ be the sub-tournament of $T$ induced on the vertices of 
   $O_{T}(v)$.  If $x\in O_{T}(v)$, then $O_{T}(v)\cap O_{T}(x)=O_{W}(x)$ 
   and so since $T$ is out-quadrangular, 
   $d_{W}^{+}(x)=|O_{W}(x)|=|O_{T}(v)\cap O_{T}(x)|\neq 1$.  
   \qed \end{pf*}
   
   \medskip\indent
   Applying Theorem~\ref{odeg23} to the dual of a tournament we get 
   the following theorem.
   
   \begin{theorem}\label{ideg23} Let $T$ be an in-quadrangular 
   tournament and choose $v\in V(T)$.  Let $W$ be the sub-tournament 
   of $T$ induced on $I(v)$.  Then $W$ contains no vertices of 
   in-degree $1$. \end{theorem}

   The only tournaments on $2$ or $3$ vertices are the single arc, 
   the $3-$cycle and the transitive triple, each of which contain a 
   vertex of out-degree $1$ and a vertex of in-degree $1$.  Therefore, 
   Theorems~\ref{odeg23} and \ref{ideg23} give us the following three 
   corollaries.

   \begin{corollary}\label{tourn8}  If $T$ is an out-quadrangular 
   tournament with 
   $\delta^{+}(T)\geq 2,$ then \\ $\delta^{+}(T)\geq 4.$  
   \end{corollary}

   \begin{corollary}\label{tourn8a}  If $T$ is an in-quadrangular 
   tournament with $\delta^{-}(T)\geq 2$, then \\ $\delta^{-}(T)\geq 4.$ 
   \end{corollary}

   \begin{corollary}\label{tourn8b} If $T$ is a quadrangular 
   tournament with $\delta^{+}(T)\geq 2$ and \\ $\delta^{-}(T)\geq 2$, 
   then $\delta^{+}(T)\geq 4$ and $\delta^{-}(T)\geq 4.$ 
   \end{corollary}

   \section{Quadrangularity in regular tournaments}
   
   In this section we look at regular tournaments, and how this 
   requirement affects quadrangularity.  We will see that regularity 
   actually makes the job of determining whether or not a tournament 
   is quadrangular a bit easier.  We also restate the problem of 
   whether or not a rotational tournament is quadrangular in a 
   more number theoretic context.  First we need the following 
   definitions and proposition.  Let $D$ be a digraph, and $x\in 
   V(D)$.  The \em closed outset \em of $x$, denoted $O[x]$ is the 
   set $O(x)\cup\{x\}$.  Similarly, the \em closed inset \em of $x$ 
   is $I[x]=I(x)\cup\{x\}$.

    \begin{proposition}  Let $T$ be a tournament on $n$ vertices, then 
    $T$ is \\ in-quadrangular
    if and only if for all $u\neq v\in V(T),$  $|O[u]\cup 
    O[v]|\neq n-1.$ \end{proposition}
    
    \begin{pf*}{Proof.} Note that since $T$ is a tournament 
    $I(x)=V(T)-O[x]$ for all $x\in V(T)$.  Since
    $T$ is in-quadrangular if and only if $|I(u)\cap 
    I(v)|\neq 1$ for all $u\neq v\in V(T),$ we have that $T$ is 
    in-quadrangular if and only if for all $u\neq v\in V(T)$ 
    \begin{eqnarray*} 1 & \neq & |I(u)\cap I(v)| \\
	& = & |(V(T)-O[u])\cap (V(T)-O[v])| \\
	& = & |V(T)-(O[u]\cup O[v])| \\
	& = & n - |O[u]\cup O[v]|. 
	\end{eqnarray*}
    Thus, $T$ is in-quadrangular if and only if  $|O[u]\cup O[v]|\neq 
    n-1$ for all $u\neq v\in V(T).$  \qed \end{pf*}

    \medskip\indent
    From this proposition, we can see that a
    tournament $T$ is  quadrangular if and only if for all 
    $u\neq v\in V(T)$, $|O(u)\cap O(v)|\neq 1$ and $|O[u]\cup O[v]|\neq 
    n-1$.  
    
    \begin{theorem}\label{regquada} A regular tournament is quadrangular 
    if and only  
    if it is out-quadrangular or in-quadrangular. \end{theorem}
    
    \begin{pf*}{Proof.} Let $T$ be a regular tournament on 
    $n=2k+1$ vertices.  Note 
    that for any two distinct vertices $x$ and $y$ in $T$, $|O[x]\cap 
    O[y]|=|O(x)\cap O(y)|+1 $ since either $x\beats y$ or $y\beats 
    x.$  Thus for any two $x\neq y\in V(T),$ 
    \begin{eqnarray*}
	|O[x]\cup O[y]| & = & 2k+2-|O[x]\cap O[y]|\\
	& = & n+1 -|O[x]\cap O[y]|\\
	& = & n+1 - 1 -|O(x)\cap O(y)|\\
	& = & n -|O(x)\cap O(y)|.
    \end{eqnarray*}
    Therefore, $|O[x]\cup O[y]|=n-1$ if and only if $|O(x)\cap 
    O(y)|=1.$  Thus, $T$ is out-quadrangular if and only if it is 
    in-quadrangular, and so it is quadrangular if and only if 
    out-quadrangular or in-quadrangular.
    \qed \end{pf*}

    \medskip\indent
    The following results give us a sufficient condition for a 
    regular tournament to be out-quadrangular in terms of the domination 
    number of the tournament.

    \begin{theorem}\label{regquad}  
	If $T$ is a regular tournament with 
    $\gamma(T)\geq 4,$ then $T$ is out-quadrangular. \end{theorem} 
    
    \begin{pf*}{Proof.} Let $T$ be a regular tournament on $2k+1$ vertices 
    with 
    $\gamma(T)\geq 4$.
    Assume to the 
    contrary that, $T$ is not out-quadrangular.  Then there exist 
    $u\neq v\in 
    V(T)$ such that  $|O(u)\cap O(v)|=1$ 
    Let $w\in V(T)$ be the single vertex in $O(u)\cap O(v)$,
    and without loss of generality assume that $u\beats v.$  
    So $|O(u)\cup O(v)|=2k-1$, since $O(u)\cap O(v)=\{w\}$.  So,$ 
    |O[u]\cup 
    O[v]|=2k$ since $u\beats v$.  Thus, there is only one 
    vertex in $T-\{u,v\}$ which is not dominated by $u$ or $v$, call 
    it $x$.  Then every vertex in $T$ is either one of $u,v,x$ or 
    dominated by one of $u,v,x$, hence $\{u,v,x\}$ form a dominating 
    set of order $3$ in $T$.  This contradicts our assumption that 
    $\gamma(T)\geq 4$. Thus, $T$ is 
    out-quadrangular.  \qed \end{pf*}
    
    \medskip\indent
    From Theorem~\ref{regquada}, we get the following corollary.
    
    \begin{corollary} If $T$ is a regular tournament with 
    $\gamma(T)\geq 4,$ then $T$ is quadrangular. \end{corollary}

   \medskip\indent
   Let $S$ be a set of $k$ integers between $1$ and $2k$ 
   such that if $i,j\in S,$ $i+j\neq 
   2k+1.$  Construct a digraph on $2k+1$ vertices labeled 
   $0,1,\ldots,2k,$ with $i\beats j$ if and 
   only if $j-i\pmod{2k+1}\in S.$  By our choice of $S,$ if 
   $i\beats j$ then $j$ does not beat $i$.  Also, 
   this digraph must have $(2k+1)k = \binom{2k+1}{2}$ arcs.  Thus, 
   this is a tournament.  Such a tournament is called a \em rotational 
   tournament, \em and the set $S$ is called its \em symbol. \em  We 
   denote by $U_{n}$ the rotational tournament whose symbol is 
   $\{1,2,\ldots,\frac{n-1}{2}\}$.  In \cite{Lundgren}, Fisher, 
   Lundgren, Merz and Reid show that if a tournament on $n$ vertices 
   has an $n$-cycle as its domination graph, then it is isomorphic 
   to $U_{n}$.
    
    \begin{theorem}\label{notU}  If $T$ is a rotational tournament on 
    $n=2k+1$ 
    vertices, then either 
    $T$ is isomorphic to $U_{n}$ or for all $u\neq v\in V(T)$, $O(u)\cap 
    O(v)\neq \emptyset$. \end{theorem}
    
    \begin{pf*}{Proof.} Let $T$ be a rotational tournament
    which is not isomorphic 
    to $U_{n}.$  Suppose, to the contrary, that there exist vertices 
    $u\neq v\in V(T)$ such that $O(u)\cap O(v)=\emptyset.$  
    Without loss of generality suppose $u=0$.  Now, 
    since $T$ is regular and $O(u)\cap O(v)=\emptyset,$ $|O[u]\cup 
    O[v]|=2k+2 - 1 - 0=2k+1=n.$  Thus, $u$ and $v$ form a dominant 
    pair.  Since 
    $T$ is rotational, $O(u+i)\cap O(v+i)=\emptyset$ for all $i$. 
    This says that $\{i,i+v\pmod{n}\}$ forms a 
    dominant pair for all $i.$  This means that the domination graph 
    of $T$ is a cycle, and the only tournaments with this property 
    are isomorphic to $U_{n}.$  \qed \end{pf*}

    \medskip\indent
    Pick $n>3$.  Then for $0,\frac{n-3}{2}\in 
    V(U_{n})$, $|O(0)\cap O(\frac{n-3}{2})|=1$.  So, $U_{n}$ 
    is not quadrangular for any $n>3$, and we get the following 
    corollary.

    \begin{corollary}\label{prop1} If $T$ is rotational and quadrangular 
    with 
    $|V(T)|>3$, then \\ $O(u)\cap O(v)\neq \emptyset$ for all 
    $u\neq v\in V(T).$  \end{corollary}

    \begin{theorem} Let $T$ be a rotational tournament on $n>3$ 
    vertices, with symbol $S$.  Then, $T$ is quadrangular if and only 
    if for all integers $m$ with $1\leq m \leq \frac{n-1}{2}$ 
    there exist distinct subsets $\{i,j\},\{k,l\}\subseteq S$ such 
    that \\ 
    $(i-j)\equiv(k-l)\equiv m\pmod{n}.$ \end{theorem}

    \begin{pf*}{Proof.} Pick $u\neq 0\in V(T),$ and suppose $S$ 
    has the property  
    stated in the theorem. We show that $|O(u)\cap 
    O(0)|\neq 1$.  For $x\in V(T)$ we can use the 
    rotational property of $T$ to map $x$ to $0$.  Then, as $u$ was 
    arbitrarily chosen, for any vertex $y\neq x$ we have $|O(x)\cap 
    O(y)|=|O(0)\cap O(u)|\neq 1$.  So $T$ will be out-quadrangular 
    and hence quadrangular since $T$ is regular.  If $u\leq 
    \frac{n-1}{2}$ there exist sets $\{i,j\},\{k,l\}\subseteq S$ such 
    that $(i-j)\equiv(k-l)\equiv u\pmod{n}.$  So, $i-u\equiv j\pmod{n}$ 
    and $k-u\equiv l\pmod{n}$  Thus, $j,l\in O(u).$  Further, $j,l\in 
    O(0)$ since $j,l\in S.$  Note, $j\neq l$ for otherwise, $i=k$ 
    contradicting $\{i,j\}$ and $\{k,l\}$ being distinct sets.  
    Thus, $|O(u)\cap O(0)|\geq 2.$  If $u\geq 
    \frac{n-1}{2}$ then $-u\leq\frac{n-1}{2}$ and so there exist 
    sets $\{i,j\},\{k,l\}\subseteq S$ such 
    that $(i-j)\equiv(k-l)\equiv -u\pmod{n}.$  So,
    $(j-i)\equiv(l-k)\equiv u\pmod{n},$ and the argument is the same.
    
    \medskip\indent
    Now, assume that $T$ is quadrangular.  Then by 
    Corollary~\ref{prop1}, $|O(u)\cap O(v)|\geq 1$ for all $u,v\in 
    V(T).$  Thus, for all $u,v\in V(T),$ we must have that $|O(u)\cap 
    O(v)|\geq 2.$  In particular, for all $m\in V(T),$ we must have 
    that $|O(0)\cap O(m)|\geq 2.$  Since $O(0)=S,$ there must be 
    at least $2$ elements of $S$ say $j,l$ such that $j,l\in O(m).$  
    So, there must exist $i,k\in S$ such that $i-m\equiv j\pmod{n}$ and 
    $k-m\equiv l\pmod{n}.$  Note this makes $\{i,j\}$ and $\{k,l\}$ 
    the sets stated in the theorem.  Also, if $m\geq \frac{n-1}{2}$, 
    then   
    $-m\pmod{n}\leq \frac{n-1}{2}$ and the argument in the previous 
    paragraph shows that the sets which work for $-m$ 
    also work for $m.$  
    \qed \end{pf*}

    \medskip\indent
    This theorem lets us restate the existence question for 
    quadrangular rotational tournaments $(n>3)$ as the following: \\
    For which odd integers $n$ does there exist a set 
    of size $\frac{n-1}{2}$ such that if $i\in S,$ $-i\pmod{n}\not\in S$ 
    and for all integers $1\leq m \leq \frac{n-1}{2},$ there exist 
    distinct sets $\{i,j\},\{k,l\}\subseteq S$ such that 
    $(i-j)\equiv (k-l)\equiv m \pmod{n}$?  
    
    \smallskip\indent
    The smallest such $n$ is $11$ with $S=\{1,3,4,5,9\}.$  
    In fact, one can generalize this set 
    and verify that for $n\equiv 3\pmod{4}$ the set 
    $$S=\left\{i:1\leq i \leq n-2, i \mbox{ is odd }, i\neq 
    \left(\frac{n+3}{2}\right)\right\}\cup 
    \left\{\left(\frac{n-3}{2}\right)\right\}$$
    is the symbol for a quadrangular rotational tournament.

    \section{Further work and open problems}
    
    There is still work to be done in this area.  Though we have 
    given a number of necessary conditions and some classifications 
    for a tournament to be quadrangular we lack  
    constructions of quadrangular tournaments.  Also, a number of our 
    classifications required a lower bound on the domination number 
    of a sub-tournament.  Quite a bit of work has been done on 
    domination in tournaments already
    (for example \cite{Lundgren,Fisher}).
    However, finding tournaments with a given 
    domination number is still an open problem.

    \medskip\indent
    Quadrangularity is a nice property for examining the 
    structure of orthogonal matrices, but not all quadrangular 
    digraphs are the digraph of an orthogonal 
    matrix.  For instance, as mentioned above, the rotational 
    tournament with symbol $S=\{1,3,4,5,9\}$ is a quadrangular 
    tournament.  
    However, the adjacency matrix of this tournament is the 
    incidence matrix of the $(11,5,2)$ design which, as is shown in 
    \cite{Klee}, cannot be the pattern of a real orthogonal matrix.   
    So, stronger necessary conditions should be explored.

    \medskip\indent
    In a coming paper, we address some of these issues and 
    determine for which orders quadrangular tournaments exist.   
    We also explore a more restrictive necessary condition known 
    as strong quadrangularity.


\begin{thebibliography}{1}

\bibitem{Klee}
L.~B. Beasley, R.~A. Brualdi, B.~L. Shader, Combinatorial orthogonality, in:
  R.~A. Brualdi, S.~Friedland, V.~Klee (Eds.), Combinatorial and
  Graph-Theoretical Problems in Linear Algebra, Vol.~50 of The IMA Volumes in
  Mathematics and its Applications, Springer-Verlag, New York, 1993, pp.
  207--218.

\bibitem{Gibson}
P.~M. Gibson, G.-H. Zhang, Combinatorially orthogonal matrices and related
  graphs, Linear Algebra Appl. 282 (1998) 83--95.

\bibitem{Reid/Thomassen}
K.~B. Reid, C.~Thomassen, Edge sets contained in circuits, Israel J. Math 24
  (1976) 305--319.

\bibitem{Kus}
K.~Zyczkowski, M.~Kus, W.~Slomczynski, H.-J. Sommers, Random unistochastic
  matrices, J. Phys. A: Math. Gen.  (28 March 2003) 3425--3450.

\bibitem{Lundgren}
D.~C. Fisher, J.~R. Lundgren, S.~K. Merz, K.~B. Reid, The domination and
  competition graphs of a tournament, J. Graph Theory 29 (1998) 103--110.

\bibitem{Fisher}
D.~C. Fisher, J.~R. Lundgren, S.~K. Merz, K.~B. Reid, Domination graphs of
  tournaments and digraphs, Congr. Numer. 108 (1995) 97--107.

\end{thebibliography}
    \end{document}